\documentclass[12]{article}
\usepackage{fleqn}
\usepackage{graphicx}
\begin{document}
\Large
\begin{center}
{\bf Geometric Hyperplanes: Desargues Encodes Doily}
\end{center}
\vspace*{.2cm}
\begin{center}
Metod Saniga
\end{center}
\vspace*{.0cm} \normalsize
\begin{center}
$^{1}$Astronomical Institute, Slovak Academy of Sciences,
SK-05960 Tatransk\' a Lomnica\\ Slovak Republic\\
(msaniga@astro.sk)

\end{center}

\vspace*{.0cm} \noindent \hrulefill

\vspace*{.1cm} \noindent {\bf Abstract}

\noindent It is shown that the structure of the generalized quadrangle of order two is fully encoded in the properties of the Desargues configuration. A point of the quadrangle is represented by a geometric hyperplane of the  Desargues configuration and its line by a set of three hyperplanes such that one of them is the complement of the symmetric difference of the remaining two and they all share a pair of non-collinear points.

\vspace*{-.1cm} \noindent \hrulefill

\vspace*{.3cm}

Certain small $v_3$-configurations exhibit a lot of remarkable properties and seem to be intricately linked with each other. A special relation in this respect is that where one configuration fully determines or, so to speak, encodes the properties of the other. A particularly nice illustration of such a link is the Fano plane (i.\,e., the unique $7_3$-configuration) encoding the structure of the split Cayley hexagon of order two \cite{psm}.  The purpose of this research note is to point out --- also in view of possible physical applications --- a similar relation between other two well-known $v_3$-configurations, the Desargues configuration and the generalized quadrangle of order two.

The Desargues configuration is one of the most important point-line incidence structures. It consists of 10 points and 10 lines, with three points on a line and three lines through a point. There exist altogether ten such $10_3$-configurations (see, for example, \cite{gru}). The Desargues configuration is, unlike the others, flag-transitive and the only one where for {\it each} of its points the three points that are not collinear with it lie on a line. This configuration is depicted in Figure 1, left, in a form showing its automorphism of order three (after Polster, \cite{pol}). The generalized quadrangle of order two, GQ$(2,\,2)$, is another prominent point-line geometry. It has 15 points/lines, also with three points/lines on a line/per point. And although there are as many as 245\,342 such $15_3$-configurations, GQ$(2,\,2)$, apart from being flag-transitive, is the only one that is {\it triangle-free} (see, for example, \cite{bbp,bgpz}). Its well-known picture, showing an automorphism of order five and often dubbed the doily, is given in Figure 1, right.

\begin{figure}[pth!]
\centerline{\includegraphics[width=4.cm,height=4.cm,clip=]{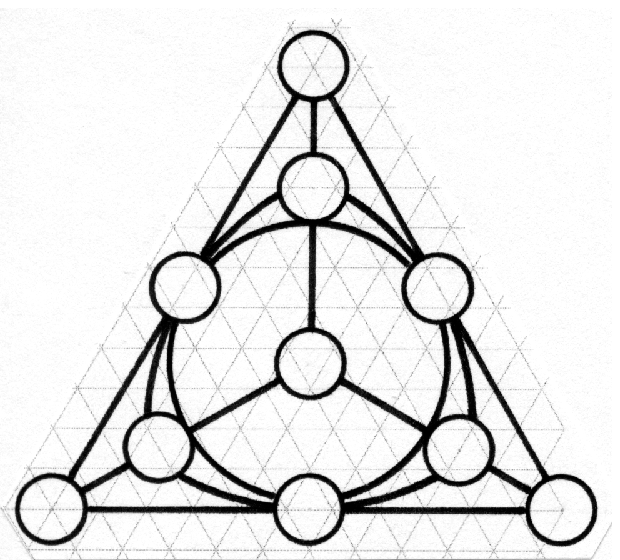}
\includegraphics[width=4.cm,height=4.cm,clip=]{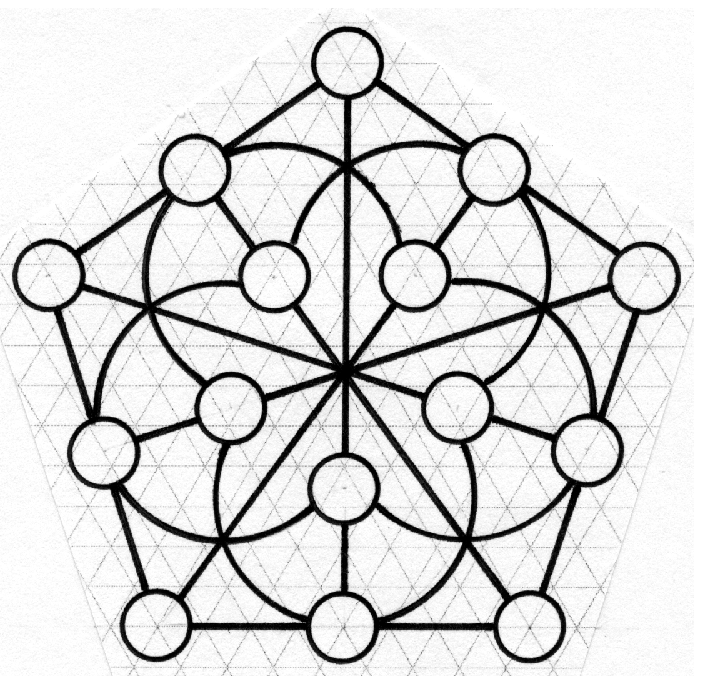}}
\vspace*{.2cm} \caption{{\it Left}: --- A picture of the Desargues configuration: small circles stand for  its points, whereas its lines are represented by triples of points on common straight segments (6), arcs of circles (3) and a big circle. {\it Right}:
--- A picture of the generalized quadrangle of order two, the doily: as in the preceding case, circles stand for its points, while its lines are given by  point triples joined by the same segments (10) and/or arcs of circles (5).}
\end{figure}

\begin{figure}[pth!]
\centerline{\includegraphics[width=8.cm,clip=]{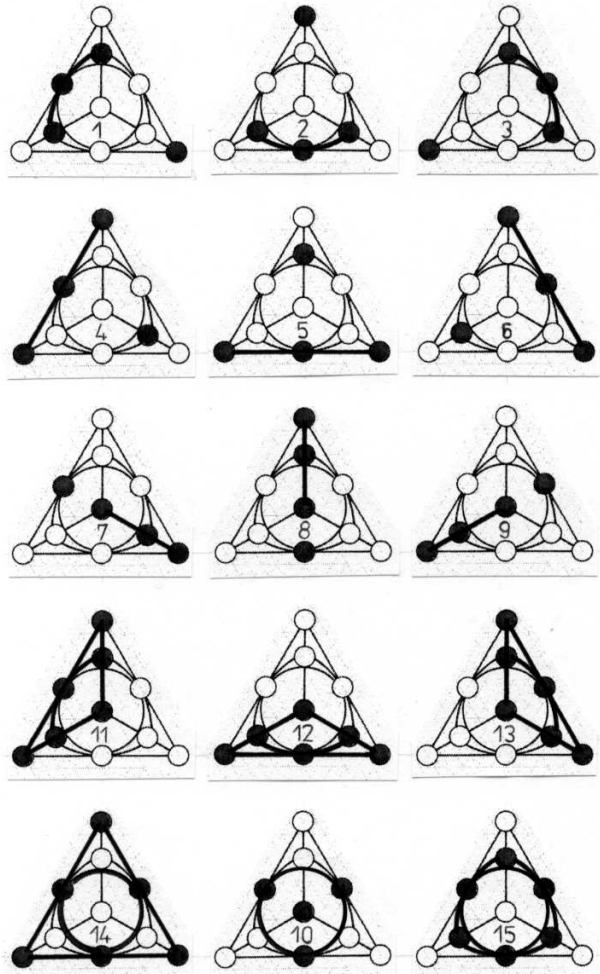}}
\vspace*{.2cm}
\caption{The 15 geometric hyperplanes (filled circles) of the Desargues configuration.}
\end{figure}

\begin{figure}[pth!]
\centerline{\includegraphics[width=14.cm,clip=]{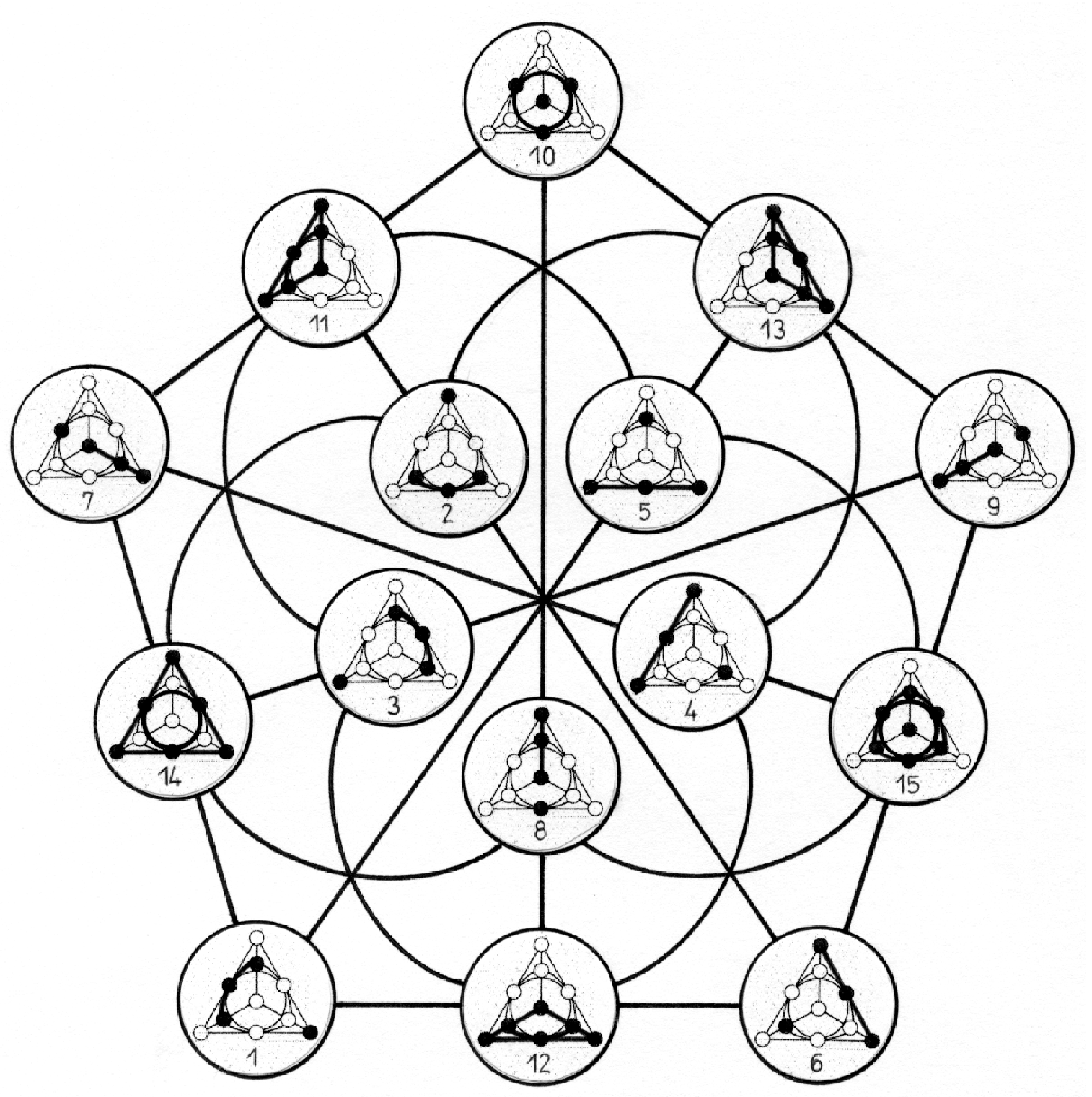}}
\vspace*{.2cm}
\caption{A diagrammatical illustration of the representation of GQ$(2,\,2)$ in terms of geometric hyperplanes of the Desargues configuration. (The reader familiar with the structure of the doily will recognize that, e.\,g., the hyperplanes 11 to 15 form an ovoid.)}
\end{figure}

These two configurations can be linked in a number of ways. The one we shall briefly describe here, believing it is novel, relies on the concept of a geometric hyperplane \cite{ron}. A geometric hyperplane of a point-line incidence geometry is a point-subset of the geometry such that each line of the geometry meets the subset in one or all points. Employing Figure 1, left, it is quite straightforward to verify that the Desargues configuration is endowed with 15 distinct geometric hyperplanes. They form two different classes and are sketched in Figure 2; any of the first ten (Nos. 1 to 10) comprises four points, namely a point and the three points non-collinear with it, whilst each of the remaining five hyperplanes (Nos. 11 to 15) features six points on three lines forming a triangle and a fourth line cutting the three in triangle's ``midpoints."  Let us now construct the point-line incidence geometry where a point is a geometric hyperplane of the Desargues configuration and a line comprises three hyperplanes such that any of them is the complement of the symmetric difference of the remaining two and they all have a pair of {\it non}-collinear points in common; a brief inspection of Figure 3 confirms that this geometry is indeed isomorphic to GQ$(2,\,2)$.
\\ \\
{\bf Acknowledgement}\\
This work was partially supported by the VEGA grant agency projects 2/0092/09 and 2/0098/10.

\end{document}